\documentclass{amsart}
\usepackage{graphicx}
\usepackage{graphicx,esint}
\usepackage{amssymb,amscd,amsthm,amsxtra}
\usepackage{latexsym}
\usepackage{epsfig}
\usepackage{color} 
\newtheorem{thm}{Theorem}[section]

\newtheorem{prop}[thm]{Proposition}
\theoremstyle{definition}
\newtheorem{defn}[thm]{Definition}
\theoremstyle{remark}
\newtheorem{rem}[thm]{\bf Remark}
\numberwithin{equation}{section}

\newcommand{\R}{\mathbb R}
\newcommand{\be}{\begin{equation}}
\newcommand{\ee}{\end{equation}}

\newcommand{\ep}{\eps}

\newcommand{\eps}{\varepsilon}

\newcommand{\comment}[1]{}

\begin{document}

\title[Two-phase free boundary problems with nonstandard growth]{Two-phase free boundary problems for operators with nonstandard growth}

\author[Fausto Ferrari]{Fausto Ferrari}
\address{Dipartimento di Matematica dell'Universit\`a di Bologna, Piazza di Porta S. Donato, 5, 40126 Bologna, Italy.}
\email{\tt fausto.ferrari@unibo.it}
\author[Monica Jacob]{Monica Jacob}
\address{Departamento  de
Ma\-te\-m\'a\-ti\-ca, Facultad de Ciencias Exactas y Naturales,
Universidad de Buenos Aires, (1428) Buenos Aires, Argentina.}
\email{\tt mljacob167@gmail.com}
\author[Claudia Lederman]{Claudia Lederman}
\address{IMAS - CONICET and Departamento  de
Ma\-te\-m\'a\-ti\-ca, Facultad de Ciencias Exactas y Naturales,
Universidad de Buenos Aires, (1428) Buenos Aires, Argentina.}
\email{\tt  clederma@dm.uba.ar}
\thanks{This work was partially
supported by a grant from the IMU-CDC and Simons Foundation}
\thanks{F. F. was partially supported by
 PRIN 2022 7HX33Z - CUP J53D23003610006, Pattern formation in nonlinear phenomena and  INDAM-GNAMPA 2024 project: {\it Free boundary
problems in noncommutative structures and degenerate operators} CUP E53C23001670001.  }
\thanks{C. L. was partially supported by the grants CONICET PIP 11220220100476CO 2023-2025,  UBACYT 20020220200056BA and ANPCyT PICT 2019-00985. C. L. wishes to thank the Department of Mathematics of the University of Bologna, Italy, for the kind hospitality.
}
\dedicatory{To our dear friend Sandro Salsa, a free boundaries master, on the  occasion of his 75th birthday}

\begin{abstract}
In this paper we focus the attention on free boundary problems ruled by partial differential equations with nonstandard growth, presenting in particular some recent results. The interest in these problems stems  from the diverse  applications that motivate their study and from the challenging
mathematical difficulties they pose.
\end{abstract}

\maketitle

\section{Introduction}\label{chap:1}

In the modern theory of partial differential equations, the scientific interest is often addressed towards  nonlinear problems---those for which a linear combination of solutions does not turn out to be a solution. 
These equations arise in the modeling of numerous disciplines,  including engineering, physics, chemistry, and biology.

 The classical theory assumes that the  nonlinear terms involved in the differential equation are of power-law type. Nevertheless, this kind  of equations 
 sometimes is insufficient to describe certain complex phenomena. Hence a more general class is being employed for the modeling of inhomogeneous materials (with different properties at different points), anisotropic materials (with different behavior in different directions), diffusive processes in strongly inhomogeneous porous media, and in models 
ranging from modern materials science to meteorology and image reconstruction.

 This type of nonlinear equations, where the growth law is more general than a power-law, is known in literature as {\it nonstandard growth equations}---a term  that has been introduced in \cite{Ma1} and \cite{Ma2}. It is a broad class of equations that
 requires a variety of tools in order to construct a suitable functional analytical framework (see  \cite{CGSW, DHHR, HH, MR}).

The aim of this paper is focusing the attention on  free boundary problems ruled by partial differential equations with nonstandard growth. We will concentrate on a class of PDEs with nonstandard growth known as {\it partial differential equations with variable exponents}. They have 
a wide range of applications, such as the modeling of non-Newtonian fluids (e.g., electrorheological \cite{R} or thermorheological \cite{AR}). Other application fields include nonlinear elasticity \cite{Z1}, image restoration \cite{AMS, CLR}, the modeling of electrical conductors \cite{Z2}, as well as the filtration process of gases in inhomogeneous porous media \cite{AS}, to cite a few.

The natural functional spaces for studying the type of equations considered here are the Lebesgue and Sobolev spaces  with variable exponents $L^{p(\cdot)}(\Omega)$ and  $W^{1,p(\cdot)}(\Omega)$, where $1< p(x) <\infty$, which generalize the classical Lebesgue and Sobolev spaces $L^{p}(\Omega)$ and $W^{1,p}(\Omega)$. We refer to \cite{DHHR} for a detailed study of these spaces. Here, we simply present  some basic tools (see Section \ref{bas-tool} and Appendix \ref{appA1}).

We start our exposition by introducing some models involving nonlinear partial differential equations with variable exponents arising from different applications.

We  then address some recent results for one and two-phase free boundary problems for operators with nonstandard growth, 
focusing our attention on free boundary problems for a particular family of PDEs with nonstandard growth of the form
\begin{equation}\label{LW3.1}
\mathrm{div}A(x,u,\nabla u) = B(x,u,\nabla u)\quad \mbox{in}\,\Omega.
\end{equation}
This family of equations is given by $A(x, s, \eta)$ which satisfies\begin{equation}\label{acotac}
\lambda_0|\eta|^{p(x)-2}|\xi|^2 \leq \sum_{i,j=1}^n \frac{\partial A_i}{\partial \eta_j}(x, s, \eta)\xi_i\xi_j \leq \Lambda_0|\eta|^{p(x)-2}|\xi|^2, \quad \xi \in \mathbb{R}^n,    
\end{equation}
where $0 < \lambda _0 \leq \Lambda _0$, and has a right-hand side given by $B(x, s, \eta) \not\equiv 0$.
In the particular case where $p(x) \equiv 2$ the equation is uniformly elliptic since it satisfies\begin{equation*} 
\lambda_0 |\xi|^2 \leq \sum_{i,j=1}^n \frac{\partial A_i}{\partial \eta_j}(x, s, \eta)\xi_i\xi_j \leq \Lambda_0 |\xi|^2, \quad \xi \in \mathbb{R}^n.    
\end{equation*}

However, $p(x)$ is variable and satisfies $1<p_{\min}\le p(x) \le p_{\max}<\infty$. Consequently, equation \eqref{LW3.1} is not uniformly elliptic. Indeed, \eqref{acotac} indicates that it is singular in the regions where $1 < p(x) < 2$ \, \big ($\lambda_0|\eta|^{p(x) - 2} \to \infty$ when $|\eta| \to 0$\big ), and degenerate in those where $p(x) > 2$ \, \big ($\Lambda_0|\eta|^{p(x) - 2} \to {0}$ when $|\eta| \to {0}$\big ).

The family under consideration was dealt with in the  work of X. Fan \cite{Fan}, where the $C^{1,\alpha}$ regularity  of bounded weak solutions was proved. Then, in \cite{Wo}, it was obtained a Harnack Inequality for this family.

More recently, in  \cite{LW5}, under the assumptions  in \cite{Fan}, several results for these equations were proved. Among them, existence of  solution, comparison, uniqueness, maximum principle and boundedness of solutions.

The method used to obtain the existence of solution consists of showing that, given $\phi\in W^{1,p(\cdot)}(\Omega)$, there exists a minimizer $u$ of the energy functional\begin{equation}\label{LW3.2}
\mathcal{J}{_\Omega}(v)=\int_\Omega{F(x,v,\nabla v)\,dx,}    
\end{equation}
in $\phi+W^{1,p(\cdot)}_0(\Omega)$, where $F$ is an appropriate  function. Consequently, $u$ is a solution to \eqref{LW3.1}, with
\begin{align*}\label{LW3.2'}
A(x,s,\eta)&=\nabla_\eta F(x,s,\eta),\\
B(x,s,\eta)&=F_s(x,s,\eta).
\end{align*}

We refer to  \cite{LW5} for a detailed study of this family. Further novel results for this family can be found in \cite{FJLM}.

A prototype operator in the \eqref{LW3.1}-\eqref{acotac} class is given by
\begin{equation*}\label{LW1.1'}
\Delta_{p(x)}u: =\mathrm{div}(|\nabla u(x)|^{p(x)-2}\nabla u)= f(x).
\end{equation*}
The operator $\Delta_{p(x)}$, called the $p(x)$-Laplacian, extends the Laplacian, which is obtained in the particular case where $p(x)\equiv2$, and the $p$-Laplacian, in the case where $p(x)\equiv p>1$. 

For additional examples in this class, we refer to \cite{LW5} and \cite{FJLM}.

We here introduce some one and  two-phase free boundary problems for nonstandard growth operators arising in different contexts. In particular we  present recent regularity results for viscosity solutions to these problems and their free boundaries, proven in a series of   works (see \cite{FL1, FL2, FL3,  FL4, FL5, FLS}).

The paper is organized as follows. In Section \ref{chap:2} we introduce some models involving PDEs with variable exponents. Sections \ref{sect-bas-fb} and \ref{bas-tool} deal with some basics on free boundary problems and on operators with variable exponents.
In Sections \ref{sect-one-ph} and \ref{sec-two-ph} we  present free boundary regularity results for viscosity solutions to one and two-phase free boundary problems for nonstandard growth operators. In Section \ref{lip-reg} we discuss optimal regularity results for 
viscosity solutions to two-phase free boundary problems. Finally, in  Appendix \ref{appA1},  we
collect some  results on  Lebesgue and Sobolev spaces with variable exponent.

\section{Some models with variable exponents}\label{chap:2}
The  systematic study of differential equations with variable exponents was motivated by the description of different models such as electrorheological and thermorheological fluids, image processing, or robotics. For an introduction to these topics, we refer to \cite{HHLN} and \cite{RaRe}. In this section we present some important models involving this type of equations.


\subsection{Image processing}
 Chen, Levine, Rao  dealt in \cite{CLR} with an application to image reconstruction. In fact, consider an input $\textit{I}$---the received image---which corresponds to shades of gray in a domain $\Omega\subset \mathbb{R}^2$. Suppose that $\textit{I}$ consists of the true image corrupted by noise. Assume moreover that the noise is additive. That is, $\textit{I} = \textit{T} + \eta$, where $\textit{T}$ is the true image and $\eta$ is the noise. The effect of noise can be removed by smoothing the input. Mathematically, smoothing corresponds to minimizing the energy
\[
E_1(v) = \int_{\Omega} \big(|\nabla v(x)|^2 + (v(x) - I(x))^2\big) \, dx.
\]
However, smoothing destroys fine image details, so this method is unsatisfactory. A better approach is the total variational smoothing. Since an edge in the image results in a very large gradient, the level sets around the edge are very dissimilar, so this method does a good job by preserving edges. More precisely, total variational smoothing corresponds to minimizing the energy
\[
E_2(v) = \int_{\Omega} \big(|\nabla v(x)| + (v(x) - I(x))^2\big) \, dx.
\]
This type of smoothing  preserves edges but, unluckily, also creates edges where none existed in the original image. Hence, by observing $E_1$ and $E_2$, Chen, Levine, and Rao suggested that an appropriate energy would be
\[
E(v) = \int_{\Omega}\Big (\frac{|\nabla v(x)|^{p(x)}}{p(x)}  + \frac{1}{2}(v(x) - I(x))^2\Big ) \, dx,
\]
where $1 \leq p(x) \leq 2$. This function should be close to 1 where there might possibly  be edges, and close to 2 where there might possible  be no edges. The approximate location of the edges can be determined  by simply smoothing the input data and observing where the gradient is large.

It is not hard to see that a minimizer of $E$ is a solution to the variable exponent equation $$\Delta_{p(x)}u = u-I.$$

\subsection{Non-Newtonian fluids}\label{f no new}
\subsubsection{\bf Introduction}\label{Intro fnn}
The equations that model a homogeneous, incompressible fluid (see, for instance, \cite{GFA, Lions}) are
\begin{equation*} \label{Intro fnn-b}
\begin{cases}
 u_t + (u \cdot \nabla)u - \mathrm{div}(S(x,Du)) +\nabla \pi &= f, \\
 \mathrm{div}u &=0,
\end{cases}     
\end{equation*}
where $u$ is the velocity, $\pi$ is the pressure, and $f$ is the external force.

The first equation describes the conservation of momentum, and the second is the condition of incompressibility.

The stress tensor is defined as
\[
T=-\pi I + S(x,Du),
\]
where
\[
Du= \frac{1}{2}(\nabla u + {\nabla u}^t)
\]
(i.e., the symmetric part of the gradient of $u$) is called the deformation tensor.

In the particular case in which  $S(x,Du)$ has the form \[
S(x,Du)= \nu Du,
\]
with $\nu$ a positive constant, the fluid is called Newtonian. An example of this type of fluids is water.

If $S(x,Du)$ depends on $Du$ in a nonlinear manner, the fluid is called non-Newtonian. Examples of this type of fluids are paint, blood, ketchup, and toothpaste.

\subsubsection{\bf Electrorheological fluids}
Electrheological fluids (we refer to Ružička, \cite{R} for a monograph on this topic) are fluids characterized by their ability to undergo significant changes in their mechanical properties when an electric field is applied. This property can be exploited in technological applications, such as clutches, shock absorbers, and rehabilitation equipment, to name a few. Winslow \cite{W} is credited with first observing the behavior of electrorheological fluids in 1949. Early fluids suffered from major constraints, such as their abrasive nature and suspension instability, and the enormous voltage requirements necessary for a significant change in material properties. Major advances have been made to overcome these obstacles, and materials  making the aforementioned devices possible are nowadays available.

The constitutive equation for the motion of an electrorheological fluid is
\begin{equation}\label{el reo}
    u_t  + (u \cdot \nabla)u - \mathrm{div} S(x,Du) + \nabla\pi = f, 
\end{equation}
where $u$ is the fluid velocity, $\pi$ is the pressure, $f$ is the external force, and the stress tensor $T$ has the form
\[
T=-\pi I + S(x,Du),
\]
\[
S(x, Du) = \mu(x) \left[1 + |Du(x)|^2\right]^{\frac{p(x)-2}{2}} Du(x),
\]
where $Du = (\nabla u + \nabla u^t)/2$ is the symmetric part of the gradient of $u$, and $\mu$ and $p(x)$ depend on the electric field.

Let us remark that the highest order differential term in \eqref{el reo} is given by
\[ \mathrm{div} \Big( \mu(x) \left[1 + |Du(x)|^2\right]^{\frac{p(x)-2}{2}} Du(x) \Big). \]

\subsubsection{\bf Thermorheological fluids}
In recent years, various nonlinear constitutive relations for the stress tensor have been proposed. In \cite{RaRu}, Rajagopal and Ružička discussed mathematical models of electrorheological fluids where, as a growth condition, the exponent depends on the electric field.
In \cite{AR}, the authors are devoted to the analysis of steady flows of fluids that are strongly influenced by the temperature field, rather than by an external electric field---the so-called thermorheological fluids.
More precisely, a system of equations describing a stationary  thermoconvective flow of a non-Newtonian fluid is studied. It is assumed that the stress tensor $T$ has the form
\[ T = - \pi I + \big ( \mu(\theta) + \tau(\theta)|D(u)|^{p(\theta)-2}\big ) D(u), \]
where $u$ is the velocity vector, $\pi$ is the pressure, $\theta$ is the temperature, $\mu$, $p$, and $\tau$ are given coefficients that depend on the temperature, and 
\[ D(u) = \frac{1}{2}(\nabla u + \nabla u^t), \]
is the deformation tensor.

The equations considered in \cite{AR} are the following: $\Omega$ is an open and bounded subset of $\mathbb{R}^n$, $n = 2,3$ and the functions \(\theta (x)\), \(u(x) = (u_1, ..., u_n)\), \(\pi(x)\) satisfy the equations
\[ (u \cdot \nabla)b(\theta) = \Delta\theta + g, \]
\[ \theta(x) = \theta_1(x), \quad x \in \partial\Omega, \]
\[ (u \cdot \nabla)u = \mathrm{div} \Big ( \big ( \mu(\theta) + \tau(\theta)|D(u)|^{p(\theta)-2} \big) D(u) \Big ) - \nabla \pi + f, \]
\[ \mathrm{div} u = 0, \]
\[
u(x) = 0, \quad x \in \partial\Omega. \]

In these equations $u, \pi$ and $\theta$ are respectively the velocity field, the pressure and the temperature of the fluid, $D(u)$ 
is the deformation tensor and $f$ is the external force. The coefficients $b$, $\mu$, $p$ and $\tau$ depend on the temperature \(\theta\).

\subsection{\bf Non-homogeneous porous media}
Equations of the type
\begin{equation}\label{med cont}
 u_t - \mathrm{div} (|u|^{\gamma(x,t)}\nabla u) = f \quad \text{en} \quad \Omega \times [0,T],
\end{equation}
with variable exponent $\gamma(x,t)$, where $\Omega \subset \mathbb{R}^n, n=3$, appear naturally in continuum mechanics. Consider the motion of an ideal barotropic gas---that is, its density depends only on the pressure---through 
a porous medium (e.g., sand, clay). Let $\rho$ be the  gas density, $V$ the velocity, and $p$ the pressure. The motion is governed by the mass conservation law 
\[
\frac{\partial \rho}{\partial t} + \mathrm{div}(\rho V) = 0,
\]
Darcy's law, which for an inhomogeneous medium has the form
\[
V = -k(x)\nabla p,
\]
where $k(x)$ is a given matrix, and the equation of state $p = P(\rho)$. It is generally assumed that \(P(s) = \mu s^{\alpha}\) with $\mu$ and $\alpha$ constants. The above conditions then lead to the equation for the density \(\rho\)
\[
\frac{\partial \rho}{\partial t} = \frac{\mu \alpha}{1 + \alpha} \mathrm{div}(k(x) \nabla \rho^{1+\alpha}).
\]
If it is further assumed that $p$ may explicitly depend on \((x, t)\) and has the form $p = \mu \rho^{\gamma(x,t)}$, the equation for \(\rho\) becomes
\[
\frac{\partial \rho}{\partial t} = \mu\,\mathrm{div} (k(x) \rho {\nabla \rho^{\gamma(x,t)}}),
\]
and can be written in the form
\begin{equation}\label{rho sub t}
\frac{\partial \rho}{\partial t} = \mu\,\mathrm{div} \Big ( k(x) \gamma\, \rho{^\gamma} \nabla \rho + (\rho^{\gamma+1} \, \text{ln} \rho) k(x) \nabla\gamma \Big).
\end{equation}
In \cite{AS} equation \eqref{rho sub t} is derived and equation \eqref{med cont}, which is a simplified version of \eqref{rho sub t}, is studied.

\subsection{\bf Thermistors modeling} A thermistor is an electric conductor with conductivity depending on the temperature. They have many applications, for instance, as temperature sensors or  circuit breakers (see \cite{HRS, Z2} and the references therein).

In \cite{Z2}, the following thermistor model was considered. Let $\Omega\subset \mathbb{R}^n$ be a smooth bounded domain. The system 
\begin{equation}
\left\{ 
\begin{array}{ll}
-\mbox{div}(|\nabla u(x)|^{\sigma(\theta(x))-2}\nabla u(x) )=\tilde{f},\:\:\: u_{|\partial \Omega}=0
\  & \  \\ 
-\Delta \theta(x)=\lambda |\nabla u(x)|^{\sigma(\theta(x))},\:\:\:\theta_{| \partial \Omega}=0
\end{array}
\right.  \label{thermistor}
\end{equation}
gives a joint description of the electric field (with potential $u$) and the temperature $\theta$ in a thermistor.

Here $\lambda>0$ is a parameter and $\sigma:[0,\infty)\to \mathbb{R}$ is a  function such that, for some constants $a,b$, there holds that $1<a\leq\sigma(s)\leq b$ for $s\in [0,\infty).$

We now assume that the temperature $\theta$ is known and denote $p(x)=\sigma(\theta (x))$ then, the problem takes the form
\begin{equation}
\left\{ 
\begin{array}{ll}
\mbox{div}(|\nabla u(x)|^{p(x)-2}\nabla u(x) )=-\tilde{f},\:\:\:\mbox{in}\:\:\Omega
\  & \  \\ 
u=0,\:\:\:\mbox{on}\:\:\partial \Omega,\\
 |\nabla u(x)|=g(x):=\left(-\frac{\Delta \theta(x)}{\lambda}\right)^{\frac{1}{p(x)}},\:\:\:\mbox{in}\:\: \Omega.&\\
\end{array}
\right.  \label{thermistorfbp}
\end{equation}
That is, we obtain the following one-phase free boundary problem
\begin{equation}  \label{fbthermi}
\left\{
\begin{array}{ll}
\Delta_{p(x)} u = f, & \hbox{in $\Omega^+(u):= \{x \in \Omega : u(x)>0\}$}, \\
\  &  \\
|\nabla u|= g, & \hbox{on $F(u):= \partial \Omega^+(u) \cap
\Omega,$} 
\end{array}
\right.
\end{equation}
where $\Delta_{p(x)} u:= \mbox{div}(|\nabla u(x)|^{p(x)-2}\nabla u(x) )$ and $f:=-\tilde{f}.$

\section{Some basics on free boundary problems}\label{sect-bas-fb}
In Sections \ref{sect-one-ph} and \ref{sec-two-ph} we will present results on problems that, among other motivations, arise in the study of minimizers of the energy
\begin{equation}\label{energ-p(x)}
J(v)={\displaystyle\int_\Omega\Big(\frac{|\nabla
v|^{p(x)}}{p(x)}+\lambda_+\chi_{\{v>0\}}+\lambda_-\chi_{\{v\leq 0\}}+fv\Big)\,dx},
\end{equation}
 where  $1<p(x)<\infty$ and $\lambda_+>\lambda_-\ge 0$ are given numbers. When we consider minimizers of a functional like \eqref{energ-p(x)}, we are dealing with functions that are solutions to a free boundary problem in some appropriate sense. 
In particular, the simplest case is given by 
$$J_{2, 0}(v):={\displaystyle\int_\Omega\Big(\frac{|\nabla
v|^{2}}{2}+\lambda_+\chi_{\{v>0\}}\Big)\,dx},$$
where  $\lambda_-= 0,$ $p(x)\equiv 2$ and $f\equiv 0$, that is, $J= J_{2, 0}$  is the so-called Bernoulli functional. A first step in the study of this problem is proving that when $u$ is a nonnegative minimizer of $J_{2, 0}$---in a suitable set of functions on which $J_{2, 0}$ is well posed--- there holds that $u$  satisfies  $\Delta u=0$ in $\Omega^+(u):=\{x\in \Omega,\quad u>0\}$. Much more technical effort is required to obtain that the same minimizer satisfies $|\nabla u|=\sqrt{2\lambda_+}$ on $F(u):=\partial\Omega^+(u)\cap \Omega$ in some suitable sense, see \cite{AC} and \cite{ACF2}. 

In this way, we are lead to an overdetermined problem of the form
\begin{equation}  \label{fbjust}
\left\{
\begin{array}{ll}
\Delta u = 0, & \hbox{in $\Omega^+(u):= \{x \in \Omega : u(x)>0\}$}, \\
\  &  \\
|\nabla u|=\sqrt{2\lambda_+}, & \hbox{on $F(u):= \partial \Omega^+(u) \cap
\Omega.$} 
\end{array}
\right.
\end{equation}
Notice that the second condition can be written equivalently as $u_{\nu}=\sqrt{2\lambda_+}$ on $F(u)$, where $u_{\nu}$ denotes the normal derivative in the inward direction to $\Omega^+(u)$. Indeed,  \eqref{fbjust}  may be read as the Euler-Lagrange equations associated with minimizers of the functional $J_{2, 0}.$ 

Hence,  even the well-posedness of the problem requires extra care. In fact, since $F(u)$ is an unknown of the problem,  it is not clear a priori in what sense it is possible to define a solution to  \eqref{fbjust}, with no information on the  set $F(u)$. Just to fix ideas, let us observe that the function $u(x,y)=x^2-y^2$ is harmonic in the set $\{(x,y)\in \mathbb{R}^2:\:\:x>0,\:\: x^2-y^2>0\}$, but at $(0,0)$---a point belonging to the boundary of this set---it is not defined classically the normal $\nu$. More precisely, to deal with a problem like \eqref{fbjust}, first we must decide what we mean by a solution. This issue has been successfully addressed in the pioneering works of Luis Caffarelli in \cite{C1,C2,C3}, see also \cite{CS}, where the notion  of viscosity solution to free boundary problem \eqref{fbjust} was introduced. 

In fact, it is worth  recalling that a continuous function $u$ in an open set $\Omega$ is harmonic in a viscosity sense if it is both subharmonic and superharmonic in a viscosity sense. To be precise,  $u$ is subharmonic (superharmonic) in a viscosity sense when: if  $x_0\in \Omega$,  $\varphi\in C^2(\Omega)$ and $u-\varphi$ has a local maximum (minimum) at $x_0$, then $\Delta\varphi (x_0)\geq 0 \, (\leq 0).$ Concerning the free boundary condition in \eqref{fbjust}, it is asked to be verified  in the following viscosity sense. 
A point  $x_0\in F(u)$ is said to be regular from the right if it is possible to touch it with an inner ball $B(y)\subset\Omega^+(u)$,  centered at $y$ (analogously from the left). In this case, we define the inner normal at $x_0$ as $\nu:=\frac{y-x_0}{|y-x_0|}$.  So that, $u$ is a supersolution (subsolution) if the condition at the free boundary in \eqref{fbjust} is satisfied, at every  regular point from the right (left) $x_0\in F(u)$,  as follows:    $u(x)=\alpha\langle x-x_0,\nu\rangle+o(|x-x_0|)$ in $B(y)$ ($B(y)^c$) as $x\to x_0,$ with $\alpha\ge 0$ and  $\alpha^2\leq  2\lambda_+ \,  (\geq 2\lambda_+).$ 

As a consequence, it appears natural to extend the previous definition of solution to one-phase free boundary problem to many other operators for which the viscosity notion of solution makes sense. For instance, to problem
\begin{equation}  \label{fbjustfl}
\left\{
\begin{array}{ll}
\mathcal{F}(D^2 u,x) = 0, & \hbox{in $\Omega^+(u):= \{x \in \Omega : u(x)>0\}$}, \\
\  &  \\
|\nabla u|=\sqrt{2\lambda_+}, & \hbox{on $F(u):= \partial \Omega^+(u) \cap
\Omega,$} 
\end{array}
\right.
\end{equation}
where $\mathcal{F}$ is a fully nonlinear uniformly elliptic operator. Namely, if  $\Omega\subset\mathbb{R}^n$ is an open set,  $\mathcal{F}:\mathbb{S}\times\Omega\to \mathbb{R}$ is continuous and there exist $0<\lambda\leq \Lambda$ such that
$$
\lambda \|N\|\leq \mathcal{F}(M+N,x)-\mathcal{F}(M,x)\leq \Lambda\|N\|,
$$
for every $M,N\in \mathbb{S}^n:=\{P:\:\:P=P^T\:\:\mbox{is a}\:\:n\times n\:\:\mbox{real matrix}\}$ with $N\geq 0$, where
$$
\|M\|=\max\{|e|:\:\: e\:\:\mbox{is an eigenvalue of}\:\:M\},
$$
 we say that $\mathcal{F}$ is uniformly elliptic with constants $\lambda, \Lambda$, see \cite{CC}.

Analogously, the same problem can be stated for the $p$-Laplace operator, $p\in (1,\infty),$
\begin{equation}  \label{fbjustfl-b}
\left\{
\begin{array}{ll}
\Delta_pu = 0, & \hbox{in $\Omega^+(u):= \{x \in \Omega : u(x)>0\}$}, \\
\  &  \\
|\nabla u|={\left({\frac{p}{p-1}\lambda_+}\right)}^{1/p} & \hbox{on $F(u):= \partial \Omega^+(u) \cap
\Omega,$} 
\end{array}
\right.
\end{equation}
where $\Delta_pu=\mbox{div}(|\nabla u|^{p-2}\nabla u)$. We point out as well that, whenever $\varphi\in C^2$ and $\nabla \varphi\not =0$, there holds
$$
\Delta_p\varphi=|\nabla \varphi (x)|^{p-2}\left(\Delta \varphi(x)+(p-2)\Delta_{\infty}^N\varphi(x)\right),
$$ 
where $\Delta_{\infty}^N\varphi(x)=\langle D^2\varphi (x)\frac{\nabla \varphi (x)}{|\nabla\varphi(x)|},\frac{\nabla \varphi (x)}{|\nabla\varphi(x)|}\rangle.$

Hence, a notion of viscosity solution for the $p$-Laplace operator can be introduced in like manner, \cite{JJ, JLM}, and it can be applied to free boundary problems, \cite{LN1,LN2} and \cite{LR}.
In this framework,  classical research subjects concern both the optimal global regularity of solutions $u$, and the optimal regularity of their free boundaries $F(u).$ This discussion can be extended to two-phase situations (i.e., allowing functions to change sign), which introduces  additional difficulties. See \cite{CS}, \cite{V}, and \cite{FLS} for an overview containing further bibliography on this topic. 

\section{Nonstandard growth equations: some basic tools}\label{bas-tool}

In this section we present some of the basic  tools  that are useful when dealing with operators with variable exponents.

Let $p :\Omega \to  [1,\infty)$ be a measurable bounded function,
called a variable exponent on $\Omega.$ In particular, we denote:
$$p_{\max} = {\rm
ess sup} \,p(x),\:\: \mbox{and}\:\: p_{\min} = {\rm ess inf} \,p(x).$$
 
 The variable exponent Lebesgue space $L^{p(\cdot)}(\Omega)$ is defined as the set of all measurable functions $u:\Omega \to \mathbb{R}$ for which
the modular $\varrho_{p(\cdot)}(u) = \int_{\Omega} |u(x)|^{p(x)}\,
dx$ is finite.

The Luxemburg norm on this space is defined by
$$
\|u\|_{L^{p(\cdot)}(\Omega)} = \|u\|_{p(\cdot)}  = \inf\{\lambda >
0: \varrho_{p(\cdot)}(u/\lambda)\leq 1 \}.
$$
This norm makes $L^{p(\cdot)}(\Omega)$ a Banach space.
There holds the following relation between $\varrho_{p(\cdot)}(u)$
and $\|u\|_{L^{p(\cdot)}}$:
\begin{equation*}
\begin{split}
&\min\Big\{\Big(\int_{\Omega} |u|^{p(x)}\, dx\Big)
^{1/{p_{\min}}}, \Big(\int_{\Omega} |u|^{p(x)}\, dx\Big)
^{1/{p_{\max}}}\Big\}\le\|u\|_{L^{p(\cdot)}(\Omega)}\\
 &\leq  \max\Big\{\Big(\int_{\Omega} |u|^{p(x)}\, dx\Big)
^{1/{p_{\min}}}, \Big(\int_{\Omega} |u|^{p(x)}\, dx\Big)
^{1/{p_{\max}}}\Big\}.
\end{split}
\end{equation*}

Moreover, the dual of $L^{p(\cdot)}(\Omega)$ is
$L^{p'(\cdot)}(\Omega)$ with $\frac{1}{p(x)}+\frac{1}{p'(x)}=1$.

 $W^{1,p(\cdot)}(\Omega)$ denotes the space of measurable
functions $u$ such that $u$ and the distributional derivative
$\nabla u$ are in $L^{p(\cdot)}(\Omega)$. The norm
$$
\|u\|_{1,p(\cdot)}:= \|u\|_{p(\cdot)} + \| |\nabla u|
\|_{p(\cdot)}
$$
makes $W^{1,p(\cdot)}(\Omega)$ a Banach space.

The space $W_0^{1,p(\cdot)}(\Omega)$ is defined as the closure of
the $C_0^{\infty}(\Omega)$ in $W^{1,p(\cdot)}(\Omega)$.

We include in an Appendix at the end of the paper additional useful results on these spaces.

Let us point out some important differences---with respect to classical operators---found when dealing with operators with variable exponent, due to the presence of the exponent  $p(x)$.

In fact, if $u\geq 0$ is a weak solution to 
\begin{equation}\label{p(x)-homog}
\Delta_{p(\cdot)}u=0,\:\:\:\text{ in } \Omega,
\end{equation}
then, there exists a positive constant $C$---that depends on $u$---such that, for any $B_{4R}(x_0)\subset\subset \Omega$,
\begin{equation}\label{harn-homog}
\sup_{B_R(x_0)}u\leq C(\inf_{B_R(x_0)}u+ R).
\end{equation}
The dependence of $C$ on $u$ can not be removed, see \cite{HKLMP}. We remark that  inequality \eqref{harn-homog} corresponds to  Harnack inequality. However, its structure is deeply different from the standard one---still valid when $p(x)\equiv p$ in \eqref{p(x)-homog}---where $C$ does not depend on $u$, and in the right hand side there is no $R$ present. Nevertheless, this kind of nonstandard growth operators still satisfy a maximum principle, although in this case it is not a straightforward consequence of  Harnack inequality, as usually happens.

In addition, if $u\geq 0$ is a weak solution to 
$$
\Delta_{p(\cdot)}u=f,\:\:\: \text{ in }\Omega,
$$
with $f\in L^{q}(\Omega)$, then there exist constants $C>0$ and $\mu\ge 0$ such that, for any $B_{4R}(x_0)\subset\subset \Omega$,
$$
\sup_{B_R(x_0)}u\leq C(\inf_{B_R(x_0)}u+ R+\mu R).
$$
Here  $C$ depends on $u$ and $\mu$ depends on $||f||_{L^q(B_{4R}(x_0))}$ (among other dependences, which appear in a complicate but clear way), see \cite{Wo}.

{}From this rough description, it is possible to perceive some of the extra difficulties that appear when dealing with operators with nonstandard growth.

\section{One-phase free boundary problems for nonstandard growth operators}\label{sect-one-ph}

In this section we  present results obtained in \cite{FL1, FL2} (see also \cite{FLS}) on a one-phase  free boundary problem governed by the $p(x)$-Laplacian  with non-zero right hand side. More precisely, 
for a  function $p$ such that $1<p(x)<\infty$, the problem is the following:
\begin{equation}  \label{fb}
\left\{
\begin{array}{ll}
\Delta_{p(x)} u:=\mbox{div} (|\nabla u|^{p(x)-2}\nabla u)= f, & \hbox{in $\Omega^+(u):= \{x \in \Omega : u(x)>0\}$}, \\
\  &  \\
|\nabla u|= g, & \hbox{on $F(u):= \partial \Omega^+(u) \cap
\Omega.$} 
\end{array}
\right.
\end{equation}
This problem comes out from limits of a singular perturbation problem with
forcing term as in \cite{LW1}, where  solutions to  
\eqref{fb}, arising in flame propagation with nonlocal and electromagnetic 
effects, are analyzed. On the other hand, \eqref{fb} appears by minimizing the functional
\begin{equation}\label{AC-energy}
J(v)=\int_{\Omega}\left(\frac{|\nabla v|^{p(x)}}{p(x)}+\lambda(x)\chi_{\{v>0\}}+f(x)v\right)dx
\end{equation}
studied in \cite{LW3}, as well as in the seminal paper by Alt and Caffarelli \cite{AC} in the case $p(x)\equiv 2$ and $f\equiv 0.$ We refer also to \cite{LW4}, where 
\eqref{fb} appears in the study of an optimal design problem.

The works \cite{FL1, FL2}  deal with the regularity of solutions and free boundaries of viscosity solutions to \eqref{fb}. There it is  followed the strategy introduced in the important paper by De Silva \cite{D} for one-phase linear problems.  \cite{D} was further extended to two-phase uniformly elliptic problems in 
\cite{DFS1, DFS2, DFS3}. The same technique was applied to the $p$-Laplace operator ($p(x)\equiv p$ in \eqref{fb}) for the one phase case, with $p\ge 2$, in \cite{LR}. 

	Problem \eqref{fb} was originally studied for weak (variational) solutions in the linear homogeneous case in \cite{AC}, associated to \eqref{AC-energy}. \cite{AC} was generalized to the linear case with $f\not\equiv 0$ in \cite{GS, Le}. In the homogeneous case,   to a quasilinear uniformly elliptic situation \cite{ACF1}, to the $p$-Laplacian \cite{DP} and to an Orlicz setting \cite{MW}. Finally, for \eqref{fb} with $1<p(x)<\infty$ and $f\not\equiv 0$, we refer to \cite{LW2}.

In \cite{FL2} it was  shown that viscosity solutions are locally Lipschitz continuous, which is the optimal regularity. Then it was proved that flat and Lipschitz free boundaries of  viscosity solutions are $C^{1,\alpha}$ (see \cite{FL1, FL2}).

The study of these issues brought challenging difficulties and  carrying out the strategy in \cite{D}  for the inhomogeneous $p(x)$-Laplacian     required new tools. This is a nonlinear operator that appears naturally in divergence form from minimization problems, i.e., in the form ${\rm div}A(x,\nabla u)=f(x)$, satisfying \eqref{acotac}.
It is singular in the regions where $1<p(x)<2$  and  degenerate in the ones where $p(x)>2$.

Among the novelties in these works,  some results that are crucial for Theorem \ref{Lipmain} below were proved. In particular Lemma 5.1 in \cite{FL2}, where  some  lemmas were revisited, for the case of $p_0$-harmonic functions, that are well known in the linear setting. This result concerns  first order expansions at one side regular boundary points of positive Lipschitz functions, vanishing at the boundary of a domain, and the proof applies to a general class of fully nonlinear degenerate elliptic operators.

It is worth remarking that  the arguments in  \cite{D} are based on  Harnack inequality, but this inequality for the $p(x)$-Laplacian is different from the standard one (see Section \ref{bas-tool}). In order to  apply this  strategy,  a new Harnack inequality for the inhomogeneous $p(x)$-Laplacian, for small perturbation settings, was proved (Theorem 3.2 in \cite{FL2}).

\subsection{Assumptions}\label{assump}
We let $\Omega\subset\R^n$  a bounded domain.   We
assume that $p\in C^1(\Omega)$, $1<p_{\min}\le p(x)\le p_{\max}<\infty$ and $\nabla p \in L^{\infty}(\Omega)$,
for some positive constants $p_{\min}$ and $p_{\max}$. We  also assume that $f\in C(\Omega)\cap L^{\infty}(\Omega)$ and that $g\in C^{0, \beta}(\Omega)\cap L^{\infty}(\Omega)$ with $g(x)\ge\gamma_0>0$,
for some positive constants $0<\beta<1$ and $\gamma_0$.

\subsection{Basic definitions, notation and preliminaries}
Let $u:\Omega\subset \mathbb{R}^n\to \mathbb{R}$ be a  continuous function.  We denote
\begin{equation*}
\Omega^+(u):= \{x \in \Omega : u(x)>0\},\qquad F(u):= \partial \Omega^+(u) \cap \Omega. 
\end{equation*} 
 In particular, if $u$ is a solution to a free boundary problem,  $\Omega^+(u)$ is called the {\it positive phase} of $u$, while  $F(u)$ is the so called {\it free boundary}. This terminology relies on the fact that, among all continuous functions, a solution $u$  is a function that solves a problem which requires that some conditions hold on the set $F(u)$ as well, even if this set is not  a priori known (see Section 
 \ref{sect-bas-fb}).

 We use the standard notion of $C$-viscosity  solution of a PDE (see e.g., \cite{CIL}).

\begin{defn}\label{defnweak} 
Assume that $1<p_{\min}\le p(x)\le p_{\max}<\infty$
with  $p(x)$ Lipschitz continuous in $\Omega$  and  $f\in L^{\infty}(\Omega)$.
We say that $u$
is a weak solution to $\Delta_{p(x)}u=f$ in $\Omega$ if $u\in W^{1,p(\cdot)}(\Omega)$ and,  for every  $\varphi \in
C_0^{\infty}(\Omega)$, there holds that
$$
-\int_{\Omega} |\nabla u(x)|^{p(x)-2}\nabla u \cdot \nabla
\varphi\, dx =\int_{\Omega} \varphi\, f(x)\, dx.
$$
\end{defn}

In \cite{FL1} the authors obtained the following useful result
\begin{thm}[Theorem 3.2 in \cite{FL1}] \label{weak-is-visc} Let $p$ and $f$ be as in Definition \ref{defnweak}. Assume moreover that $f\in C(\Omega)$ and $p\in C^1(\Omega)$. 
Let $u\in W^{1,p(\cdot)}(\Omega)\cap C(\Omega)$ be a weak solution to $\Delta_{p(x)}u=f$ in $\Omega$. Then $u$ is a viscosity solution to $\Delta_{p(x)}u=f$ in $\Omega$.
\end{thm}
\begin{defn}Given $u, \varphi \in C(\Omega)$, we say that $\varphi$
touches $u$ from below (resp. above) at $x_0 \in \Omega$ 
if $u(x_0)=
\varphi(x_0),$ and $u(x) \geq \varphi(x)$ (resp. $u(x) \leq
\varphi(x)$) in a neighborhood $O$ of $x_0$.
\end{defn}
\begin{defn}\label{defnhsol1} Let $u$ be a continuous nonnegative function in
$\Omega$. We say that $u$ is a viscosity solution to \eqref{fb} in
$\Omega$, if the following conditions are satisfied:
\begin{enumerate}
\item $ \Delta_{p(x)} u = f$ in $\Omega^+(u)$ 
in the weak sense of Definition \ref{defnweak}.
\item For every $\varphi \in C(\Omega)$, $\varphi \in C^2(\overline{\Omega^+(\varphi)})$. If $\varphi^+$ touches $u$ from below (resp.  above) at $x_0 \in F(u)$ and $\nabla \varphi(x_0)\not=0$, then 
$|\nabla \varphi(x_0)| \leq g(x_0)$ (resp. $ \geq g(x_0)$).
\end{enumerate}
\end{defn}
The following result is  a consequence of  Theorem \ref{weak-is-visc} 
\begin{thm}\label{defnhsol2} Let $u$ be a viscosity solution to \eqref{fb} in
$\Omega.$  Then the following conditions are satisfied:
\begin{enumerate}
\item $ \Delta_{p(x)} u = f$ in $\Omega^+(u)$ in the
viscosity sense.
\item For every $\varphi \in C(\Omega)$, $\varphi \in C^2(\overline{\Omega^+(\varphi)})$. If $\varphi^+$ touches $u$ from below (resp.  above) at $x_0 \in F(u)$ and $\nabla \varphi(x_0)\not=0$, then 
$|\nabla \varphi(x_0)| \leq g(x_0)$ (resp. $ \geq g(x_0)$).
\end{enumerate}
\end{thm}

\subsection{Regularity results for viscosity solutions and their free boundaries}
The following result refers to optimal regularity of viscosity solutions.
\begin{thm}[Optimal regularity, Theorem 1.1 in \cite{FL2}]
\label{Lip-contin} Let
$u$ be a viscosity solution to \eqref{fb}
in $B_1$. 
There exists a  constant $C>0$ such that
\begin{equation*}  
\|\nabla u\|_{L^{\infty}(B_{1/2})}\leq C.
\end{equation*}
\end{thm}
The proof relies on comparison with suitable  barrier functions  constructed in Lemma 4.2 in \cite{FL1}.

The following achievement about the fact that flat free boundaries are $C^{1,\alpha}$ was obtained
\begin{thm}[Flatness implies $C^{1,\protect\alpha}$, Theorem 1.1 in \cite{FL1}]
\label{flatmain1} Let
$u$ be a viscosity solution to \eqref{fb}
in $B_1$. Assume that  $0\in F(u),$ $g(0)=1$ and $p(0)=p_0.$   
There exists a universal constant $\bar{\varepsilon}>0$ such that, if the graph of $u$ is $\bar{\varepsilon}-$flat in $B_1,$  in the direction $e_n,$ that is
\begin{equation*}  \label{trapped}
(x_n-\bar{\varepsilon})^+\leq u(x)\leq (x_n+\bar{\varepsilon})^+, \quad x\in B_1,
\end{equation*}
\begin{equation*}  \label{pflat}
\|\nabla p\|_{L^{\infty}(B_1)}\leq \bar{\varepsilon},\quad \|f\|_{L^{\infty}(B_1)}\leq \bar{\varepsilon}, \quad [g]_{C^{0,\beta}(B_1)}\leq \bar{\varepsilon},
\end{equation*}
 then $F(u)$ is $C^{1,\alpha}$ in $
B_{1/2}$.
\end{thm}
The proof is based on an
improvement of flatness, obtained via a compactness argument which linearizes the
problem into a limiting one.  The key tool is a Harnack inequality that allows the  passage to
the limit.


Moreover, it is also possible to prove that Lipschitz free boundaries are indeed $C^{1,\alpha}$.
\begin{thm}[Lipschitz implies $C^{1,\alpha}$, Theorem 1.2 in \cite{FL2}]
\label{Lipmain} Let $u$ be a viscosity solution to \eqref{fb} in $B_1$, with
$0\in F(u)$. 
 If $F(u)$ is a Lipschitz graph in a neighborhood of $0$, then 
$F(u)$ is $C^{1,\alpha}$ in a (smaller) neighborhood of $0$.
\end{thm}
The proof follows from Theorem \ref{flatmain1} and the main
result in  \cite{LN1}, via a blow-up argument. The application
of this result required nontrivial arguments  due to the different notion of solution employed in \cite{LN1}.

\subsection{Some applications and conclusions}
We discuss shortly some applications of the results  in \cite{FL1} and \cite{FL2} and we draw some conclusions
on them. These applications correspond to three different minimization problems studied in \cite{LW1}, \cite{LW3} and \cite{LW4}.  We refer to them for motivation and related literature.
The first application is
\begin{prop}\label{applic-AC} Let $\Omega$, $p$ and $f$ be as above. Let $0<{\lambda_{\min}}\le\lambda(x)\le{\lambda_{\max}}<\infty$ with 
$\lambda\in C^{0, \beta}(\Omega)$.
Let $u\in W^{1,p(\cdot)}(\Omega)\cap L^{\infty}(\Omega)$   be a nonnegative local minimizer of 
$J(v)={\displaystyle\int_\Omega\Big(\frac{|\nabla
v|^{p(x)}}{p(x)}+\lambda(x)\chi_{\{v>0\}}+fv\Big)\,dx}$
in $\Omega$. Then,  $u$ is a viscosity solution to \eqref{fb} in $\Omega$ with $g(x)=(\frac{p(x)}{p(x)-1}\,\lambda(x))^{1/p(x)}$.
Let $x_0\in F(u)$ be such that $F(u)$ is a Lipschitz graph in a neighborhood of $x_0$, then  $F(u)$ is $C^{1,\alpha}$ in a neighborhood of $x_0$.
Let  $x_0\in F(u)$ be such that $F(u)$ has a normal in the measure theoretic sense, then $F(u)$ is $C^{1,\alpha}$ in a neighborhood of $x_0$.
Moreover, there is a subset $\mathcal{R}$  of $F(u)$ which is locally a
$C^{1,\alpha}$ surface. $\mathcal{R}$ is open and
dense in $F(u)$ and the remainder has $(n-1)-$dimensional Hausdorff measure zero.
\end{prop}
\begin{rem} \label{applic-optim}
In \cite{LW1}, \cite{LW3} and \cite{LW4} minimization and  optimization problems were considered, for  energies associated to the $p(x)$-Laplacian,  and  the same conclusions as in Proposition \ref{applic-AC} is obtained for them from the application of \cite{FL1} and \cite{FL2}.
\end{rem}
\begin{rem} \label{rem-concl-applic}
In Proposition \ref{applic-AC} and Remark \ref{applic-optim}, the $C^{1,\alpha}$ regularity results on $F(u)$ under the Lipschitz assumption follow from Theorem \ref{Lipmain} and are new.
The rest the $C^{1,\alpha}$ regularity results on $F(u)$ in Proposition \ref{applic-AC} and Remark \ref{applic-optim}, which follow from Theorem \ref{flatmain1}, were already obtained in \cite{LW3} and \cite{LW4}, from the application of the results in \cite{LW2}, but under different assumptions on $f$ and $p$.

In fact, the results in \cite{FL1}---inspired by \cite{D}---require that $f\in C(\Omega)\cap L^{\infty}(\Omega)$ and $p\in C^1(\Omega)$ and Lipschitz, whereas the results in \cite{LW2}---inspired by  \cite{AC}---require that $f\in L^{\infty}(\Omega)\cap W^{1,q}(\Omega)$ and $p\in W^{1,\infty}(\Omega)\cap W^{2,q}(\Omega)$, for $q>\max\{1, n/2\}$.
The reason for this difference relies on the fact that in the approach in \cite{D} for viscosity solutions the estimates are obtained by comparison with
suitable barriers. In the approach in \cite{AC} for variational solutions, estimates on $|\nabla u|$ close to the free boundary are obtained by looking for an equation for $v=|\nabla u|$, which requires to differentiate the equation and thus, stronger assumptions on the data $f$ and $p$ are needed in \cite{LW2}.
\end{rem}

\section{Two-phase free boundary problems for nonstandard growth operators}\label{sec-two-ph}

In this section we  consider   the following  two-phase free boundary problem:
\begin{equation}  \label{fb2}
\left\{
\begin{array}{ll}
\Delta_{p(x)} u = f, &  \hbox{in $\Omega^+(u)\cup \Omega^-(u),$} \\
\  &  \\
u_\nu^+=G(u_\nu^-,x), & \hbox{on $F(u):= \partial \Omega^+(u) \cap
\Omega,$} 
\end{array}
\right.
\end{equation}
where $\Omega  \subset \mathbb{R}^n$ is a bounded domain and
$$\Omega^+(u):= \{x \in \Omega : u(x)>0\}, \qquad\quad \Omega^-(u):=\{x \in \Omega : u(x)\leq 0\}^{\circ },$$
 while $u_\nu^+$ and $u_\nu^-$ denote the normal derivatives in the inward
direction to $\Omega^+(u)$ and $\Omega^-(u)$ respectively. 

This problem comes out  naturally from limits of singular perturbation problems with
forcing term as in \cite{LW} and  \cite{LW1}, where solutions to  
\eqref{fb2}, arising in the study of flame propagation with nonlocal and electromagnetic 
effects, are analyzed. On the other hand, nonnegative solutions to \eqref{fb2} appear in \cite{LW4} where 
an optimal design problem is studied. Problem \eqref{fb2} is also obtained 
by minimizing the following functional 
\begin{equation*}
{J}(v)=\int_{\Omega}\left(\frac{|\nabla v|^{p(x)}}{p(x)}+q(x)(\lambda_+\chi_{\{v>0\}}+\lambda_-\chi_{\{v\leq 0\}})+f(x)v\right)dx,
\end{equation*}
where $\lambda_+>\lambda_-\ge 0$ are given numbers and $q$ is a strictly positive given function. For nonnegative minimizers  we refer to \cite{LW3}. In case of minimizers without sign restriction, the two-phase problem  \eqref{fb2} is obtained with free boundary 
condition given by
\begin{equation*}
(u_\nu^+)^{p(x)}-(u_\nu^-)^{p(x)}=q(x)p(x)\frac{(\lambda_+-\lambda_-)}{p(x)-1},
\end{equation*}
under suitable assumptions (see Appendix A in \cite{FL4}, and \cite{DF})

We will present free boundary regularity results for viscosity solutions to problem \eqref{fb2} obtained in \cite{FL4}.

In fact, in \cite{FL4} it is   followed the strategy developed in \cite{DFS1, DFS2} for two-phase uniformly elliptic settings---inspired by \cite{D}---and it is proved that flat free boundaries of two-phase viscosity solutions to
\eqref{fb2} are $C^{1,\gamma}.$

Let us mention that the two-phase problem \eqref{fb2}, in the linear homogeneous case, governed by the Laplacian, i.e., when $p(x)\equiv 2$ and $f\equiv 0$,  was settled in the classical works by
Caffarelli \cite{C1,C2}. 
These results have been widely generalized to different
classes of homogeneous elliptic problems. See for example \cite{CFS, FS1,
FS2} for linear operators, \cite{AF,F1, F2, Fe1,W1, W2} for fully nonlinear
operators and \cite{LN1},\cite{LN2} for the $p$-Laplacian. The general strategy followed by these papers, however, seems not so suitable when a non-zero right hand side is 
present, as it is the  case of \eqref{fb2}.

	We also point out that, in the one-phase case, problem \eqref{fb2} with non-zero right hand side was dealt with in \cite{LW2}. There,  the $C^{1,\alpha}$ regularity
	of the free boundary near flat free boundary points was obtained, for weak (variational) solutions, following the approach in \cite{AC}. However, it is not clear how to adapt these techniques to the two-phase case.

One of the key differences between the  situation in \cite{FL4} and the one in \cite{DFS1} and \cite{DFS2} is that in these later works, given $u$ a viscosity solution to the free boundary problem, the functions $v=u-\alpha x_n$ and
$v=u-\beta x_n$ are solutions in $\Omega^+(u)\cup \Omega^-(u)$ to the same equation as $u$  which, of course, is still uniformly elliptic. This fact is repeatedly used throughout the proofs.
In contrast, in case of  problem \eqref{fb2}, such functions $v$ are viscosity solutions in $\mathcal{D}:=\Omega^+(u)\cup \Omega^-(u)$  to 
{\it an inhomogeneous equation with nonstandard growth of 
general type} of the form 
\begin{equation}\label{op-general}
{\rm div}A(x,\nabla v)=f(x)\quad \text{ in }\mathcal{D},
\end{equation}
 where $A:\mathcal{D}\times\mathbb{R}^n\to\mathbb{R}^n$ satisfies the following structure conditions:

For some positive constants $C_1,C_2,C_3, C_4$, and for  every $x\in\mathcal{D}$ and $\xi\in \mathbb{R}^n$, 
\begin{equation}\label{str-above}
|A(x,\xi)|\le C_1 |\xi|^{p(x)-1} +C_2
\end{equation}
and
\begin{equation}\label{str-below}
\langle A(x,\xi),\xi\rangle\ge  C_3|\xi|^{p(x)}-C_4.
\end{equation}
We stress that the treatment  of an equation of singular/degenerate type satisfying \eqref{op-general}, \eqref{str-above} and \eqref{str-below} is highly nontrivial, in particular when---as in the present  case---the right hand side in \eqref{op-general} is not zero.

\subsection{Assumptions}\label{assump2}
The assumptions on the function $p(x)$ are: $p\in C^1(\Omega)$, $1<p_{\min}\le p(x)\le p_{\max}<\infty$, $\nabla p \in L^{\infty}(\Omega)$,
for some positive constants $p_{\min}$ and $p_{\max}$. Moreover, it is also assumed that  $f\in L^{\infty}(\Omega)$ and $f$  continuous in 
$\Omega^+(u)\cup \Omega^-(u)$.

It is also required that the function $G:[0,\infty)\times\Omega\rightarrow(0,\infty)$ satisfies, for  $0<\hat{\beta}<L$,
\begin{itemize}
\item[(H1)]
$G(\eta,\cdot)\in C^{0,\bar\gamma}(\Omega)$ uniformly in $\eta\in [\frac{\hat{\beta}}{2},4L]$; $G(\cdot,x)\in C^{1,\bar\gamma}([\frac{\hat{\beta}}{2},4L])$ for every $x\in \Omega$ and $G\in L^\infty((\frac{\hat{\beta}}{2},4L)\times \Omega)$.
\item[(H2)]
$G'(\cdot, x)>0$ in $[\frac{\hat{\beta}}{2},4L]$ for $x\in \Omega$ and, for some $\gamma_0$  constant, $G\geq\gamma_0>0$ in $[\frac{\hat{\beta}}{2},4L]\times \Omega$.
\end{itemize}
These assumptions are complemented with  additional structural conditions  (see \cite{FL4}). This allows the inclusion of some interesting free boundary conditions, that are discussed in detail in Section 7 in \cite{FL4}.

\subsection{Basic definitions, notation and preliminaries}\label{bas-defin-2}

 For any continuous function $u:\Omega\subset \mathbb{R}^n\to \mathbb{R}$ we denote
\begin{equation*}
\Omega^+(u):= \{x \in \Omega : u(x)>0\},\quad \Omega^-(u):=\{x \in \Omega : u(x)\leq 0\}^{\circ }
\end{equation*} 
and
\begin{equation*}
F(u):= \partial \Omega^+(u) \cap \Omega. 
\end{equation*} 
We refer to the set $F(u)$ as the {\it free boundary} of $u$, while $\Omega^+(u)$ is its {\it positive phase} (or {\it side}) and $\Omega^-(u)$ is the {\it nonpositive phase}.

\begin{defn}
\label{defnhsolcv1} Let $u$ be a continuous function in $\Omega$. We say that
$u$ is a viscosity solution to \eqref{fb2} in $\Omega$, if the following
conditions are satisfied:
\begin{enumerate}
\item $\Delta_{p(x)} u = f$ in $\Omega^+(u) \cup
\Omega^-(u)$ in the weak sense of Definition \ref{defnweak}.
\item  Let $x_0 \in F(u)$ and $v \in C^2(\overline{B^+(v)}) \cap C^2(\overline{B^-(v)})$ ($B=B_\delta(x_0)$) with $F(v) \in C^2$. If $v$ touches $u$ by below (resp. above) at $x_0 \in F(v)$, then 
$$v_\nu^+(x_0) \leq G(v_\nu^-(x_0), x_0) \quad (\text{resp. $ \geq$)}.$$\end{enumerate}
\end{defn}

Next theorem follows as a consequence of  Theorem \ref{weak-is-visc}. 

\begin{thm}
\label{defnhsolcv2} Let $u$ be a viscosity solution to \eqref{fb2} in $\Omega$. Then the following
conditions are satisfied:
\begin{enumerate}
\item $\Delta_{p(x)} u = f$ in $\Omega^+(u) \cup
\Omega^-(u)$ in the viscosity sense.
\item  Let $x_0 \in F(u)$ and $v \in C^2(\overline{B^+(v)}) \cap C^2(\overline{B^-(v)})$ ($B=B_\delta(x_0)$) with $F(v) \in C^2$. If $v$ touches $u$ by below (resp. above) at $x_0 \in F(v)$, then 
$$v_\nu^+(x_0) \leq G(v_\nu^-(x_0), x_0) \quad (\text{resp. }  \geq).$$\end{enumerate}
\end{thm}

\subsection{Free boundary regularity results for two-phase viscosity solutions}\label{fb-reg-2}

Let $U_\beta$ be the two-plane solution to \eqref{fb2} when $p(x)\equiv p_0$ and $f \equiv 0 $, i.e.,
\begin{equation*}
U_\beta(x) = \alpha x_n^+ - \beta x_n^-, \quad \beta \geq 0, \quad \alpha
=G(\beta, 0).
\end{equation*}

One of the results  obtained in \cite{FL4} is the following theorem, which is the two-phase counterpart of  Theorem 1.1 in \cite{FL1} in the one-phase setting.
The main result tells us that from the flatness  of the free boundary descends the $C^{1,\protect\gamma}$ regularity of the free boundary.
\begin{thm}[Flatness implies $C^{1,\gamma}$, Theorem 1.4 in \cite{FL4}]\label{main_new-nondeg-DS-general} Let $u$ be a viscosity solution to  
\eqref{fb2} in $B_1$. Let $0<\hat{\beta}<L$. There exists a universal
constant $\bar \varepsilon>0$ such that, if 
\begin{equation*} U_{\beta}(x_n-\bar\varepsilon) \leq u(x) \leq U_\beta(x_n+ \bar\varepsilon)   \ \text{in} \ B_1\quad \text{for some }\,
0 <\hat{\beta}\leq \beta \leq L,\end{equation*} 
 and $$\|\nabla p\|_{L^\infty(B_1)} \leq \bar \varepsilon, \qquad \|f\|_{L^\infty(B_1)} \leq \bar \varepsilon,$$
\begin{equation*} [G(\eta,\cdot)]_{C^{0,\bar \gamma}(B_1)} \leq \bar \varepsilon
, \quad \text{for all }\, 0<\frac{\hat{\beta}}{2} \leq \eta \leq 4L,
\end{equation*}
then $F(u)$ is $C^{1,\gamma}$ in $B_{1/2}$. Here $\gamma$ is universal and the  $C^{1,\gamma}$ norm of $F(u)$ is bounded by a universal constant.
\end{thm}

As a consequence of Theorem \ref{main_new-nondeg-DS-general},  the following result is obtained

\begin{thm}[Flatness implies $C^{1,\gamma}$, Theorem 1.3 in \cite{FL4}]\label{main_new-nondeg-general} Let $u$ be a viscosity solution to  \eqref{fb2} in $B_1$. Let $0<\hat{\beta}<L$.  There exists a universal
constant $\bar \ep>0$ such that, if \begin{equation*}\|u - U_{\beta}\|_{L^{\infty}(B_{1})} \leq \bar \eps\quad 
\text{for some }\, 0 <\hat{\beta}\leq \beta \leq L,\end{equation*} 
 and $$\|\nabla p\|_{L^\infty(B_1)} \leq \bar \ep, \qquad \|f\|_{L^\infty(B_1)} \leq \bar \ep,$$
$$[G(\eta,\cdot)]_{C^{0,\bar \gamma}(B_1)} \leq \bar \ep
, \quad \text{for all }\, 0<\frac{\hat{\beta}}{2} \leq \eta \leq 4L,$$
then $F(u)$ is $C^{1,\gamma}$ in $B_{1/2}$. Here $\gamma$ is universal and the  $C^{1,\gamma}$ norm of $F(u)$ is bounded by a universal constant.
\end{thm}

The proof of Theorem \ref{main_new-nondeg-DS-general} is based on an
improvement of flatness, obtained via a compactness argument which linearizes the
problem into a limiting one.  The key tool is a geometric Harnack inequality
that localizes the free boundary well, and allows the rigorous passage to
the limit.

\subsection{Some comments and conclusions}\label{conc-2}

It is worth remarking that no assumption on the Lipschitz continuity of solutions is made. These regularity results are the first ones in literature for two-phase free boundary problems for the $p(x)$-Laplacian and also for two-phase problems for singular/degenerate 
operators with non-zero right hand side. They are new even when $p(x)\equiv p$, i.e., for the  $p$-Laplacian. The fact that these results hold for merely viscosity solutions allows a wide applicability.

We  would like to stress that in order to carry out this study it was necessary to overcome challenging difficulties due to its nonlinear singular/degenerate nature (see \cite{FL1, FL2, FL3, FL4, FLS}).

\section{Lipschitz continuity of two-phase viscosity solution}\label{lip-reg}

In this section we  present regularity results obtained in \cite{FL5}  for viscosity solutions to two-phase free boundary problems for the $p$-Laplacian with non-zero right hand side, where $p\in (1,\infty)$. 

In order to deal with the regularity issue for solutions to two-phase free boundary problems, it is worth starting from the simplest case, i.e., when $p= 2$.  In \cite{ACF2}, via the celebrated monotonicity formula, it was proven---in the homogeneous case---that a minimizer $u$, in a suitable class of functions, of
$$J(v)={\displaystyle\int_\Omega\Big(\frac{|\nabla
v|^{2}}{2}+\lambda_+\chi_{\{v>0\}}+\lambda_-\chi_{\{v\leq 0\}}\Big)\,dx}$$ 
is locally Lipschitz continuous. Here  $\lambda_+>\lambda_-\ge0$ are given constants.

In fact, it was shown that a minimizer $u$ satisfies: $u\in C(\Omega)$,
\begin{equation}  \label{fbtrue2}
\left\{
\begin{array}{ll}
\Delta  u = 0, & \hbox{in $\Omega^+(u)\cup \Omega^-(u)$}, \\
\  &  \\
u_\nu^+=\left({2(\lambda_+ -\lambda_-)+(u_\nu^-)^2}\right)^{1/2} & \hbox{on $F(u):= \partial \Omega^+(u) \cap
\Omega,$} 
\end{array}
\right.
\end{equation}
in a generalized sense. Then, the local Lipschitz continuity of $u$  follows from the  fact that, for every $x_0\in F(u)$, there exists  $r_0>0$ such that the function $\Phi:(0,r_0]\to \mathbb{R}$, defined as
$$
\Phi(r):=\frac{1}{r^4}\int_{B_r(x_0)}\frac{|\nabla u^+(x)|}{|x-x_0|^{n-2}}dx\int_{B_r(x_0)}\frac{|\nabla u^-(x)|}{|x-x_0|^{n-2}}dx,
$$
is bounded and monotone increasing. 

No monotonicity formula is known for the analogous problem ruled by the $p$-Laplace operator when $p\not=2$ (see \cite{DK2} for the status of the play). Nevertheless, a Lipschitz  regularity result was obtained in \cite{DK1} for minimizers of the homogeneous $p$-Bernoulli two-phase functional
\begin{equation*}
{J}_p(v)=\int_{\Omega}\left(\frac{|\nabla v|^{p}}{p}+(\lambda_+\chi_{\{v>0\}}+\lambda_-\chi_{\{v\leq 0\}})\right)dx,
\end{equation*}
for $p\in (1,\infty)$ and $\lambda_+>\lambda_-\ge0$  given constants.

On the other hand, continuing the research  in \cite{FL1}, \cite{FL2} and \cite{FL4} on  inhomogeneous one and two-phase free boundary problems governed by the $p(x)$-Laplace operator, the authors obtained a key tool 
 to address the regularity of viscosity solutions to two-phase free boundary problems ruled by the $p$-Laplace, in the inhomogeneous setting.
In fact, in \cite{FL5} the authors extended the approach introduced in \cite{DS1}, for homogeneous fully nonlinear uniformly elliptic operators, to the $p$-Laplace inhomogeneous case.

More precisely, we denote $$\Delta_{p} u:=\mbox{div} (|\nabla u|^{p-2}\nabla u),$$
where $p$ is a constant, $1<p<\infty$. We are interested in the regularity of viscosity solutions to the following two phase problem: 
\begin{equation}  \label{fbtrue}
\left\{
\begin{array}{ll}
\Delta_{p} u = f, & \hbox{in $\Omega^+(u)\cup \Omega^-(u)$}, \\
\  &  \\
u_\nu^+=G(u_\nu^-,x), & \hbox{on $F(u):= \partial \Omega^+(u) \cap
\Omega,$} 
\end{array}
\right.
\end{equation}
where $\Omega  \subset \mathbb{R}^n$ is a bounded domain and
$$\Omega^+(u):= \{x \in \Omega : u(x)>0\}, \qquad\quad \Omega^-(u):=\{x \in \Omega : u(x)\leq 0\}^{\circ },$$
 while $u_\nu^+$ and $u_\nu^-$ denote the normal derivatives in the inward
direction to $\Omega^+(u)$ and $\Omega^-(u)$ respectively and $F(u)$ is the \emph{free boundary} of $u$.

In fact, in \cite{FL5} the authors consider problem \eqref{fbtrue} for a 
the function 
$$G(t,x):[0,\infty)\times\Omega\rightarrow(0,\infty),$$
$G(\cdot, x)\in C^2(0,\infty)$ for $x\in \Omega$, that
satisfies:

\begin{itemize}
\item[(P1)]
$G(\cdot,x)$ is strictly increasing for $x\in \Omega$  and $G(t,x)\to \infty$ as $t\to \infty$, uniformly in $x\in \Omega$.
\end{itemize}

In the main result it is assumed that $G(t,x)$ behaves like $t$, for $t$ large, uniformly in $x$. More precisely, it is assumed that:

\begin{itemize}
\item[(P2)] $\frac{\partial G(t,x)}{\partial t}\to 1$, $\frac{\partial^2G(t,x)}{\partial t^2}=O(\frac{1}{t})$, as $t\to\infty$, uniformly in $x\in \Omega$,

\item[(P3)] $G(t,\cdot)\in C^{0,\bar\gamma}(\Omega)$ uniformly in $t$, for some $0<\bar{\gamma}<1$.
\end{itemize}

When $f\not\equiv 0$ and $p\neq 2$ it is assumed  as well that:
\begin{itemize}
\item[(P4)] $G(t,x)\equiv G(t)$, with $t^{-1}G(t)$ strictly decreasing, for $t$ large.
\end{itemize}

In particular, \cite{FL5} includes $G(t, x)=(g(x)+t^p)^{\frac{1}{p}}$, $p>1$, where $g\in C^{0,\bar{\gamma}}(\Omega)$, for some $0<\bar{\gamma}<1$, $g>0$,  which arises in several applications (see for instance, Appendix A in \cite{FL4}).

The assumptions on $f$ are
\begin{equation}\label{assump-f}
f\in L^{\infty}(\Omega),\qquad f\text{ is continuous in }\Omega^+(u)\cup \Omega^-(u).
\end{equation}
Nevertheless, the arguments also apply when  $f$ is merely bounded measurable, but \eqref{assump-f} is required  to avoid technicalities.

The main result proven in \cite{FL5} is the following: 

\begin{thm}[Optimal regularity, Theorem 1.1 in \cite{FL5}]\label{Lipschitz_cont}
Let $u$ be a viscosity solution to \eqref{fbtrue} in $B_1$. Assume that \eqref{assump-f} holds in $B_1$ and $G$ satisfies assumptions (P1), 
(P2) and (P3) in $B_1$. If $f\not\equiv 0$ and $p\neq 2$   assume that also (P4) holds, then 
\begin{equation*}
\|\nabla u\|_{L^\infty(B_{1/2})}\leq C (\|u\|_{L^\infty(B_{3/4})}+1),
\end{equation*}
where $C=C(n,p,\|f\|_{L^\infty(B_1)}, G)$ is a positive constant.
\end{thm}

Let us point out that this is the optimal regularity for the problem.

In \cite{FL5} the authors follow the approach introduced in  \cite{DS1} for uniformly elliptic fully nonlinear operators.

The heuristics behind the proof of Theorem \ref{Lipschitz_cont} is that, in the regime of {\it big gradients}, the free boundary condition becomes a continuity  condition for the gradient (i.e., no-jump). 
As a consequence, the desired gradient bound follows from interior $C^{1,\alpha}$ estimates for the $p$-Laplace operator and from free boundary regularity results for two-phase viscosity solutions to \eqref{fbtrue} proven in the recent paper \cite{FL4}.

These results are new even in the homogeneous situation, that is, when $f\equiv 0$. It is worth remarking that \cite{FL5} applies to merely viscosity solutions, which allows a wide applicability.

We point out that in the one  phase case, i.e., when $u\ge 0$, the local Lipschitz continuity of solutions to  \eqref{fbtrue} follows from the previous paper \cite{FL2}, where a more general one phase free boundary problem was treated (see Section \ref{sect-one-ph}).

We remark that carrying out for the inhomogeneous $p$-Laplace operator a strategy originally devised  for
 {\it uniformly elliptic operators} presents
challenging difficulties due to the type of nonlinear behavior of the $p$-Laplacian. This was already the case in the previous works \cite{FL1}, \cite{FL2} and \cite{FL4} for the treatment of this problem.  In fact, the $p$-Laplacian  is a nonlinear operator that appears naturally in divergence structure from minimization problems, i.e., in the form ${\rm div}A(\nabla u)=f(x)$, satisfying \eqref{acotac} for $p(x)\equiv p$. Dealing with this problem   in the presence of {\it a non-zero right hand side} is very delicate since,  in this case, {\it the factor $|\eta|^{p-2}$ in 
\eqref{acotac} can not be neglected}. 

In particular, the discussion at the end of Section 3.2 in \cite{DS1}, concerning the generalization of their results to some equations with zero right hand side, {\it does not apply to the inhomogeneous problem  \eqref{fbtrue}}.

The process of obtaining the local Lipschitz continuity of viscosity solutions to \eqref{fbtrue} requires to prove first the local H\"older continuity of viscosity solutions to this problem---which is interesting by itself. In fact,  this result holds under more general conditions on $G$. The precise result proven in \cite{FL5}  is the following: 
\begin{thm}[Local H\"older continuity, Theorem1.2 in \cite{FL5}]\label{holder_reg}
Let $u$ be a viscosity solution to \eqref{fbtrue} in $B_1$. Assume that, for $\sigma>0$ and $M>0$, there holds
\begin{equation}\label{G2}
\sigma t\leq G(t, x)\leq \sigma^{-1}t,\quad \mbox{for}\:\:t>M,\: x\in B_1.
\end{equation}
Then $u\in C^{0,\alpha}(B_{1/2})$, for some $0<\alpha<1$, and
$$
\|u\|_{C^{0,\alpha}(B_{1/2})}\leq C(\|u\|_{L^{\infty}(B_{3/4})}+1),
$$
where  $\alpha$ depends only $\sigma, n$ and $p$,  and $C>0$ depends only $\sigma, M, n, p$ and $\|f\|_{L^{\infty}(B_1)}$. 
\end{thm}

Let us mention that when the problem is truly nonlinear singular / degenerate with a non-zero right hand side---i.e., when $f\not\equiv 0$ and $p\neq 2$---due to the difficulties inherent to this situation, it is assumed in the  main result that (P4) also holds. We point out that assumptions of these type appear frequently in the study of free boundary problems. See, for instance, \cite{C1, C2, DFS1, DFS2}.

Let us also stress that \cite{FL5} includes the case in which $p=2$---i.e., the case of Laplace operator---relying on the work \cite{FL4}.

  Notice that, although \cite{DS1} includes Laplace operator, the case of a non-zero right hand side and a general free boundary condition $u_\nu^+=G(u_\nu^-,x)$ was not explicitly developed there. Hence, Theorem \ref{Lipschitz_cont}   gives an explicit  proof for the case $p=2$ which is different from the one that can be deduced from the monotonicity formula of \cite{CJK}---{\it only available for this particular situation}.

As already mentioned, the main tool in the approach followed in \cite{FL5}, in order to prove the local Lipschitz continuity, inspired by \cite{DS1}, is a regularity result for flat free boundaries  for viscosity solution to free boundary problem \eqref{fbtrue}. Let us point out that in the uniformly elliptic linear and fully nonlinear cases, the corresponding results were obtained in \cite{DFS1} and \cite{DFS2} respectively. It is worth mentioning that in these papers it was assumed, all along the works,  the  Lipschitz continuity of the viscosity solutions under consideration. However, in \cite{DS1} it was stated a theorem---that can be derived from the results in \cite{DFS2}---concerning the free boundary regularity for the fully nonlinear uniformly elliptic case, which holds under a nondegenerate situation, that does not assume the Lipschitz continuity of solutions.

In \cite{FL5}, the authors
make use of a free boundary regularity theorem for viscosity solutions to problem \eqref{fbtrue}, that has been proved in \cite{FL4}, where the Lipschitz continuity of solutions was not assumed. 

Let us emphasize once more that  the operator is  nonlinear singular / degenerate  with a non-zero right hand side. Thus,  challenging issues arise in \cite{FL5}, but they have been faced with the aid of the work \cite{FL4}.

\begin{rem}\label{other-problems}
We point out that obtaining a monotonicity formula seems to be more difficult when dealing with other degenerate, possibly linear, operators in noncommutative structures e.g., in Carnot groups, see \cite{FerrariForcillo1, FerrariForcillo2, FerrariForcillo3, FerrariGiovagnoli}.

Concerning the regularity of minimizers of energy functionals, other interesting related problem refers to the regularity of almost minimizers. This is a more general problem, where the functions under consideration do not satisfy a PDE. See  \cite{DavidToro2015}, \cite{DeSilvaSavin2016},  \cite{DFFV} and \cite{FerrariForcilloMerlino2025}, where the local Lipschitz continuity, which is the optimal regularity, is proved.
\end{rem}

\section*{Conflict of interest}
On behalf of all authors, the corresponding author states that there is no conflict of interest.


\appendix

\section{} \label{appA1}

\setcounter{equation}{0}

In Section \ref{bas-tool} we included some preliminaries on Lebesgue and Sobolev spaces with variable exponent. We collect here some other important results on these spaces.

\begin{thm}\label{ref}
Let $p'(x)$ such that $$\frac{1}{p(x)}+\frac{1}{p'(x)}=1.$$ Then
$L^{p'(\cdot)}(\Omega)$ is the dual of $L^{p(\cdot)}(\Omega)$.
Moreover, if $p_{\min}>1$, $L^{p(\cdot)}(\Omega)$ and
$W^{1,p(\cdot)}(\Omega)$ are reflexive.
\end{thm}

\begin{thm}\label{imb}
Let $q(x)\leq p(x)$. If $\Omega$ has finite measure, then
 $L^{p(\cdot)}(\Omega)\hookrightarrow L^{q(\cdot)}(\Omega)$
continuously.
\end{thm}

We also have the following H\"older's inequality

\begin{thm} \label{holder}
Let $p'(x)$ be as in Theorem \ref{ref}. Then there holds
$$
\int_{\Omega}|f||g|\,dx \le 2\|f\|_{p(\cdot)}\|g\|_{p'(\cdot)},
$$
for all $f\in L^{p(\cdot)}(\Omega)$ and $g\in L^{p'(\cdot)}(\Omega)$.
\end{thm}

The following version of Poincare's inequality holds

\begin{thm}\label{poinc} Let $\Omega$ be bounded. Assume that $p(x)$ is log-H\"older continuous  in $\Omega$ (that is, $p$ has a modulus of continuity $\omega(r)=C(\log \frac{1}{r})^{-1}$). For
every $u\in W_0^{1,p(\cdot)}(\Omega)$, the inequality
$$
\|u\|_{L^{p(\cdot)}(\Omega)}\leq C\|\nabla u\|_{L^{p(\cdot)}(\Omega)}
$$
holds with a constant $C$ depending only on n, $\rm{diam}(\Omega)$ and the log-H\"older modulus of continuity of $p(x)$.
\end{thm}

\smallskip

For the proof of these results and more about these spaces, see \cite{DHHR}, \cite{HH},
\cite{KR}, \cite{RaRe} and the references therein.

\bigskip



\begin{thebibliography}{9999999}

\bibitem{AMS} R. Aboulaich, D. Meskine, A. Souissi, \emph{New diffusion models in image processing}, Comput.
Math. Appl. 56 (2008), 874--882.

\bibitem{AC} H. W. Alt,
L. A. Caffarelli, \emph{Existence and regularity for a minimum problem
with free boundary}, J. Reine Angew. Math 325
(1981), 105--144.

\bibitem{ACF1}
{H.~W. Alt, L.~A. Caffarelli, A.~Friedman}, \emph{A free boundary
problem for
  quasilinear elliptic equations}, Ann. Sc. Norm. Super. Pisa Cl. Sci. (4)
  11 (1) (1984),  1--44.


\bibitem{ACF2} H.~W. Alt, L.~A. Caffarelli, A.~Friedman, \emph{Variational problems with two phases and their free boundaries}, Trans. Amer. Math. Soc. {282} (1984), no. 2, 431--461.




\bibitem{AR} S. N. Antontsev, J. F. Rodrigues, {\it On stationary thermo-rheological viscous flows}, Ann. Univ. Ferrara, Sez. VII, Sci. Mat. 
 52 (1) (2006), 19--36.


\bibitem{AS}
S. N. Antontsev, S. I. Shmarev, {\it A model porous medium equation with variable exponent of nonlinearity: existence, uniqueness and localization properties of solutions},  Nonlinear Anal. {60} (2005), 515--545.

\bibitem{AF} R. Argiolas, F. Ferrari, {\it Flat free boundaries regularity in two-phase problems for a class of fully nonlinear elliptic operators with variable coefficients}, Interfaces Free Bound. 11 (2009), no. 2, 177--199.



\bibitem{C1} L. A. Caffarelli, {\it A Harnack
inequality approach to the regularity of free boundaries. Part I:
Lipschitz free boundaries are $C^{1,\alpha}$}, Rev. Mat.
Iberoamericana 3 (1987) no. 2, 139--162.

\bibitem{C2} L. A. Caffarelli, {\it A Harnack
inequality approach to the regularity of free boundaries. Part II:
Flat free boundaries are Lipschitz}, Comm. Pure Appl. Math.
42 (1989), no. 1, 55--78.

\bibitem{C3}  L. A. Caffarelli, {\it  A Harnack inequality approach to the regularity of free boundaries. III. Existence theory, compactness, and dependence on X.} Ann. Scuola Norm. Sup. Pisa Cl. Sci. (4) 15 (1988), no. 4, 583--602.

\bibitem{CC} L. A. Caffarelli, X. Cabr\'e, {\it Fully Nonlinear Elliptic
Equations}, Colloquium Publications 43, American Mathematical
Society, Providence, RI, 1995.


\bibitem{CJK} L. A. Caffarelli, D. Jerison, C. E. Kenig, {\it Some new monotonicity theorems with applications to free boundary problems}, Ann. of Math. (2) 155 (2002), no. 2, 369--404.

\bibitem{CS} L. A. Caffarelli, S. Salsa, {\it A Geometric Approach to Free Boundary Problems.} Graduate Studies in Mathematics, 68. American Mathematical Society, Providence, RI, 2005.


\bibitem{CFS} M. C.  Cerutti, F. Ferrari, S. Salsa, {\it Two phase
problems for linear elliptic operators with variable coefficients:
Lipschitz free boundaries are $C^{1,\gamma}$}, Archive for
Rational Mechanics and Analysis, Vol 171, n.3 (2004), 329--348.




\bibitem{CLR} Y. Chen, S. Levine, M. Rao, {\it Variable exponent, linear growth functionals in image restoration},
SIAM J. Appl. Math. 66 (2006), 1383--1406.


\bibitem{CGSW}
I. Chlebicka, P. Gwiazda, A. Swierczewska-Gwiazda, A. Wroblewska-Kaminska, {\it Partial Differential Equations in Anisotropic Musielak-Orlicz Spaces}, Springer Monographs in Mathematics, Springer, 2021.


\bibitem{CIL} M. G. Crandall, H. Ishii, P. L. Lions, {\it User’s guide to viscosity solutions of second order partial differential equations}, Bull. Amer. Math. Soc. (N.S.) 27 (1) (1992),  1--67.

\bibitem{DP}{D. Danielli, A. Petrosyan}, \emph{A minimum problem with free
boundary for a
  degenerate quasilinear operator}, Calc. Var. Partial Differential Equations
  23 (1) (2005),  97--124.


\bibitem{DavidToro2015} G. David, T. Toro, {\it Regularity of almost minimizers with free boundary},
Calc. Var. Partial Differential Equations 54 (2015), 455--524.


\bibitem{D} D. De Silva, {\it Free boundary regularity for a problem with right hand side}, Interfaces and free boundaries 13 (2011), 223--238.
\bibitem{DFS1} D. De Silva, F. Ferrari, S. Salsa, {\it Two-phase problems with distributed sources: regularity of the free boundary},  Anal. PDE 7 (2014), no. 2, 267--310. 

\bibitem{DFS2} D. De Silva, F. Ferrari, S. Salsa, {\it Free boundary regularity for fully nonlinear non-homogeneous two-phase problems}, J. Math. Pures Appl. (9) 103 (2015), no. 3, 658--694.

\bibitem{DFS3} D. De Silva, F. Ferrari, S. Salsa, {\it Recent progresses on elliptic two-phase free boundary problems}, Discrete
Contin. Dyn. Syst. 39 (2019), no. 12, 6961--6978.





\bibitem{DS1} D. De Silva, O. Savin, {\it Lipschitz regularity of solutions to two-phase free boundary problems}, Int. Math. Res. Not. IMRN (2019), no. 7, 2204--2222.

\bibitem{DeSilvaSavin2016} D. De Silva, O. Savin, {\it Almost minimizers of the one-phase free boundary problem}, Comm. Partial Differential Equations 45  (2020), 913--930.


\bibitem{DHHR}
L. Diening, P. Harjulehto, P. Hasto, M. Ružička, {\it Lebesgue and Sobolev Spaces with Variable Exponents}, Lecture Notes in Mathematics 2017, Springer,  2011.

\bibitem{DFFV} S. Dipierro, F. Ferrari, N. Forcillo, E. Valdinoci, {\it Lipschitz regularity of almost minimizers in one-phase problems driven by the p-Laplace operator}, Indiana Univ. Math. J. 73 (2024), 813--854.


\bibitem{DK1} S. Dipierro, A. L. Karakhanyan,  {\it Stratification of free boundary points for a two-phase variational problem}, Adv. Math. 328 (2018), 40--81.

\bibitem{DK2} S. Dipierro, A. L. Karakhanyan, {\it A new discrete monotonicity formula with application to a two-phase free boundary problem in dimension two}, Comm. Partial Differential Equations 43 (2018), no. 7, 1073--1101

\bibitem{DF} A. Dzhugan, F. Ferrari, {\it Domain variation solutions for degenerate two phase free boundary problems},
 Mathematics in Engineering  3 (6) (2021), 1--29. 


\bibitem{Fan} X. Fan, {\it Global $C^{1,\alpha}$ regularity for variable exponent elliptic equations in divergence form}, J. Differential Equations 235 (2007), 397--417.




\bibitem{F1} M. Feldman, {\it Regularity for nonisotropic two-phase problems with Lipschitz free boundaries}, Differential Integral Equations 10 (1997), no. 6, 1171--1179.

\bibitem{F2} M. Feldman, {\it Regularity of Lipschitz free boundaries in two-phase problems for fully nonlinear elliptic equations}, Indiana Univ. Math. J. 50 (2001), no. 3, 1171--1200.


\bibitem{Fe1} F. Ferrari, {\it Two-phase problems for a class of fully nonlinear elliptic operators, Lipschitz free boundaries are $C^{1,\gamma}$},
Amer. J. Math. 128 (2006), 541--571.

\bibitem{FerrariForcillo1} F. Ferrari, N. Forcillo, {\it A new glance to the Alt-Caffarelli-Friedman monotonicity formula}, Math. Eng. 2, no. 4 (2020),  657--679.

\bibitem{FerrariForcillo2} F. Ferrari, N. Forcillo, {\it A counterexample to the monotone increasing behavior of an Alt-Caffarelli-Friedman formula in the Heisenberg group}, Atti Accad. Naz. Lincei Rend. Lincei Mat. Appl. 34 (2023), 295--306.

\bibitem{FerrariForcillo3} F. Ferrari, N. Forcillo, {\it Alt-Caffarelli-Friedman monotonicity formula and mean value properties in Carnot groups with applications}, Boll. Unione Mat. Ital. 17 (2024), 333--348.


\bibitem{FerrariForcilloMerlino2025} F. Ferrari, N. Forcillo, E. M. Merlino, {\it Regularity for almost minimizers of a one-phase Bernoulli-type functional in Carnot groups of step two}, Calc. Var. Partial Differential Equations 64 (2025), no. 4, 32 pp.


\bibitem{FerrariGiovagnoli} F. Ferrari,  D. Giovagnoli, {\it Some counterexamples to Alt-Caffarelli-Friedman monotonicity formulas in Carnot groups}, Ann. Mat. Pura Appl. (4) 204  (2025), 427--445.



\bibitem{FJLM} F. Ferrari, M. Jacob, C. Lederman, E. M. Merlino, {\it in preparation}.

\bibitem{FL1}  F. Ferrari, C. Lederman, {\it Regularity of flat free boundaries for a $p(x)$-Laplacian problem with right hand side}, 
Nonlinear Anal. 212 (2021), Article ID 112444, 25 p. 

\bibitem{FL2}  F. Ferrari, C. Lederman, {\it Regularity of Lipschitz free boundaries for a $p(x)$-Laplacian problem with right hand side}, 
J. Math. Pures Appl. 171 (2023), 26--74. 

\bibitem{FL3} F. Ferrari, C. Lederman, \emph{Free boundary regularity for a one-phase problem with non-standard growth}, Matem\'atica Aplicada, Computacional e Industrial (MACI) { 9} (2023), 211--214. 


\bibitem{FL4}  F. Ferrari, C. Lederman, {\it Regularity of flat free boundaries for two-phase $p(x)$-Laplacian problems with right hand side}, Calc. Var. Partial Differential Equations 63 (2024), no. 5, Paper No. 132, 43 pp.

\bibitem{FL5}  F. Ferrari, C. Lederman, {\it Lipschitz regularity of solutions to two-phase $p$-Laplacian free boundary problems with right hand side}, Communications in Analysis and Geometry, in press.



 \bibitem{FLS}F. Ferrari, C. Lederman, S. Salsa, {\it Recent results on nonlinear elliptic free boundary problems}, Vietnam J. Math. 50 (4)
(2022), 977--996. 

 
\bibitem{FS1} F. Ferrari, S. Salsa, {\it Regularity of the free boundary in two-phase problems for elliptic operators},
 Adv. Math. 214 (2007), 288--322.

\bibitem{FS2} F. Ferrari, S. Salsa, {\it Subsolutions of elliptic operators in divergence form and application to two-phase free boundary problems}, Bound. Value Probl. 2007, art. ID 57049, 21pp.





\bibitem{GFA}
 M. E. Gurtin,  E. Fried, L. Anand, {\it The Mechanics and Thermodynamics of Continua}, Cambridge University Press, Cambridge, 2010.

 \bibitem{GS} {B. Gustafsson, H. Shahgholian}, {\it Existence and
geometric properties of solutions of a free boundary problem in
potential theory}, J. Reine Angew. Math. 473 (1996),
137--179.


\bibitem{HH}
P. Harjulehto,  P. Hästö, {\it Orlicz Spaces and Generalized Orlicz Spaces}, Lecture Notes in Mathematics 2236, Springer,  2019. 


\bibitem{HHLN}
P. Harjulehto, P. Hästö, U. V. Le,  M. Nuortio. {\it Overview of differential equations with non-standard growth}. Nonlinear Anal. {72} (2010), 4551--4574.

\bibitem{HKLMP} P. Harjulehto, T. Kuusi, T. Lukkari, N. Marola, M. Parviainen, 
{\it Harnack's inequality for quasiminimizers with nonstandard growth conditions}.
J. Math. Anal. Appl. 344 (2008), no. 1, 504--520.

\bibitem{HRS} S. D. Howison, J. F. Rodrigues,  M. Shillor, {\it Stationary solutions to the thermistor problem}, J. Math.
Anal. Appl. 174 (2) (1993), 573--588.


\bibitem{JJ} V. Julin, P. Juutinen, {\it A new proof for the equivalence of weak and viscosity solutions for the $p$-Laplace equation}, Communications in PDE 37 (2012), no. 5, 934--946.

\bibitem{JLM} P. Juutinen, P. Lindqvist, J. Manfredi, {\it On the equivalence of viscosity solutions and weak solutions for a quasi-linear equation}, SIAM J. Math. Anal. 33 (2001), no. 3, 699--717.

\bibitem{KR}
O. Kov\'a\v{c}ik, J. R\'akosn{\'i}k, \emph{On spaces ${L}^{p(x)}$ and
  ${W}^{k,p(x)}$}, Czechoslovak Math. J {41} (1991), 592--618.



\bibitem{Le}{C.~Lederman}, \emph{A free boundary problem with a volume penalization}, Ann.
  Sc. Norm. Super. Pisa Cl. Sci. (4) 23 (2) (1996),  249--300.



\bibitem{LW} C. Lederman, N. Wolanski,  {\it A two phase elliptic singular perturbation problem with a forcing term}, J. Math. Pures Appl. (9) {86} (2006), no. 6, 552--589.


\bibitem{LW1} C.  Lederman, N. Wolanski, {\it An inhomogeneous singular perturbation problem for the $p(x)$-Laplacian}, Nonlinear Anal. 138 (2016), 300--325.

\bibitem{LW2} C.  Lederman, N. Wolanski, {\it Weak solutions and regularity of the interface in an inhomogeneous free boundary problem for the $p(x)$-Laplacian}, Interfaces Free Bound. 19 (2017), no. 2, 201--241.

\bibitem{LW3} C. Lederman, N. Wolanski, {\it Inhomogeneous minimization problems for the $p(x)$-Laplacian,} J. Math. Anal. Appl. 475 (2019), no. 1, 423--463.

\bibitem{LW4} C. Lederman, N. Wolanski, {\it An optimization problem with volume constraint for an inhomogeneous operator with nonstandard growth},  Discrete Contin. Dyn. Syst. Series A { 41 (6)} (2021), 2907--2946. 

\bibitem{LW5}
C. Lederman, N. Wolanski, {\it Lipschitz continuity of minimizers in a problem with nonstandard growth}, Contemporary PDEs between theory and modeling, Mathematics in Engineering {3 (1)} (2021), 1--17.


\bibitem{LR} R. Leit$\tilde{a}$o, G. Ricarte, 
{\it Free boundary regularity for a degenerate problem with right hand side},
Interfaces Free Bound. 20 (2018), no. 4, 577--595.


\bibitem{LN1} J. Lewis, K. Nystr{\"o}m, {\it Regularity of Lipschitz free boundaries in two phase problems for the $p$-Laplace operator}, Adv. in Math. 225 (2010), 2565--2597. 

\bibitem{LN2} J. Lewis, K. Nystr{\"o}m, {\it Regularity of flat free boundaries in two-phase problems for the $p$-Laplace operator}, Ann. Inst. H. Poincar\'e Anal. Non Lin\'aire 29 (2012), no. 1, 83--108.


\bibitem{Lions}
P. L. Lions, {\it Mathematical Topics in Fluid Mechanics. Vol. 1: Incompressible Models},
Oxford Lecture Series in Mathematics and its Applications, Oxford Clarendon Press, 1996.

\bibitem{Ma1}
P. Marcellini, {\it Regularity of minimizers of integrals of the calculus of variations with non standard growth conditions}, Arch. Ration. Mech. Anal. {105 (3)} (1989), 267--284.

\bibitem{Ma2}
P. Marcellini, {\it Regularity and existance of solutions of elliptic equations with $p, q$-growth conditions}, J. Differ. Equ. {50} (1991), 
1--30. 

\bibitem{MW}{S.~Mart\'{\i}nez, N.~Wolanski}, \emph{A minimum problem with free boundary in
  {O}rlicz spaces}, Adv. Math. 218 (6) (2008),
  1914--1971.



\bibitem{MR}
G. Mingione, V. Radulescu, {\it Recent developments in problems with nonstandard growth and nonuniform ellipticity}, J. Math. Anal. Appl. {501 (1)} (2021), 125197, 41 pp.


\bibitem{RaRe}
V. D. Radulescu, D. D. Repovs, {\it Partial Differential Equations with Variable Exponents: Variational Methods and Qualitative Analysis}, Monographs and Research Notes in Mathematics, Book 9. Chapman \& Hall / CRC Press, Boca Raton, FL, 2015.

\bibitem{RaRu}
K. P. Rajagopal, M. Ružička, {\it Mathematical modelling of electro-rheological materials}. Cont. Mech. Therm. {13} (2001), 59--78.

\bibitem{R} M. Ruzicka, {\it Electrorheological Fluids: Modeling and
Mathematical Theory}, Springer-Verlag, Berlin, 2000.




\bibitem{V} B. Velichkov, {\it Regularity of the One-phase Free Boundaries.} Lecture Notes of the Unione Matematica Italiana, 28. Springer, Cham, 2023.


\bibitem{W1}  P. Y. Wang, {\it Regularity of free boundaries of two-phase problems for fully nonlinear elliptic equations of second order. I. Lipschitz free boundaries are $C^{1,\alpha}$}, Comm. Pure Appl. Math. 53 (2000), 799--810.

\bibitem{W2}  P. Y. Wang, {\it Regularity of free boundaries of two-phase problems for fully nonlinear elliptic equations of second order. II. Flat free boundaries are Lipschitz}, Comm. Partial Differential Equations 27 (2002), 1497--1514.



\bibitem{W} 
W. M. Winslow, {\it Induced fibrillation of suspensions}, Journal of Applied Physics {20} (1949), 1137--1140.




\bibitem{Wo}
N. Wolanski, {\it Local bounds, Harnack's inequality and Hölder continuity for divergence type elliptic equations with non-standard growth}, Rev. Unión Mat. Argent. {56 (1)} (2015), 73--105.


\bibitem{Z1} V. V. Zhikov, {\it Averaging of functionals of the calculus of variations and elasticity theory}, Math. USSR. Izv. 29 (1)  (1987), 33--66.


\bibitem{Z2} V. V. Zhikov, {\it Solvability of the three-dimensional thermistor problem}, Tr. Mat. Inst. Steklova D (Differ. Uravn. i Din. Sist.) 261 (2008), 101--114.


\end{thebibliography}
\end{document}